\documentclass[twoside]{article}
\usepackage{amsfonts}
\usepackage{amsmath}
\usepackage{epsfig}
\date{}
\begin{document}

\title{\bf On the necessity of the assumptions used to prove Hardy-Littlewood and Riesz Rearrangement Inequalities}
\author{H. Hajaiej}

\maketitle

\vspace*{0.5cm}

\begin{abstract}
We prove that supermodularity is a necessary condition for the
generalized Hardy-Littlewood and Riesz rearrangement inequalities.
We also show the necessity of the monotonicity of the kernels
involved in the Riesz--type integral. 
\end{abstract}

\normalsize
\small
\renewcommand{\baselinestretch}{1.2}
\normalsize



\noindent\hrulefill

\section{Introduction}
\def\Box{\vrule height 6pt depth 1pt width 4pt}


%
%


\hspace*{1cm} Let $F: \mathbb{R} _+ ^m \to \mathbb{R} $ be a Borel
measurable function. We say that $F$ is supermodular if:
\setcounter{section}{1} \setcounter{equation}{0}
\begin{eqnarray}
\label{1.1}
 & & F(y + he_i + ke_j ) + F(y) \geq F(y + he_i ) + F(y+ ke_j
)
\\
\nonumber
 & & \quad \mbox{for every $i,j\in \{ 1, \ldots ,m\} $, $i\not= j $,
$k,h >0$,}
\end{eqnarray}
where $y= (y_1 , \ldots ,y_ m)$, and $e_i $ denotes the standard
basis vector in $\mathbb{R} ^m $, $i=1, \ldots ,m$.
\\
This condition has been first analysed in the context of integral
inequalities by G. Lorentz \cite{25}.
\\
Note that if $F$ is $C^2 $, (\ref{1.1}) is equivalent to
\begin{equation}
\label{1.2} \frac{\partial ^2 F}{\partial x_i \partial x_j} (x) \geq
0, \quad \mbox{for every $i,j\in \{ 1, \ldots ,m\} $, $i\not= j $.}
\end{equation}
The supermodularity condition has been widely used in functional analysis.
It was the main assumption to prove the generalized Hardy-Littlewood
rearrangement inequality,
\begin{equation}
\label{1.3} \int_X F(u_1 (x) , \ldots , u_m (x))\, d\mu (x) \leq
\int_X F(u_1 ^* (x) , \ldots , u_m ^* (x))\, d\mu (x),
\end{equation}
Here $u_i $ are nonnegative measurable functions that vanish at
infinity and $u_i ^* $ are their symmetric decreasing
rearrangements. The integral is taken over $X= \mathbb{R} ^n $,
$\mathbb{H}^n $ or ${\cal S} ^n $. (\ref{1.3}) has been studied by
many authors \cite{23, 25, 11, 6,
7}. Its main
applications are economics \cite{11}, chemistry \cite{30}, and
nonlinear optics where the profile of stable electromagnetic waves
traveling along a planar waveguide are given by the ground states of
the energy functional
\begin{equation}
\label{1.4}
E(u):= \frac{1}{2} \int_{\mathbb{R}} \left( u^{' 2 } -
 G(|x|,u) \right) \, dx ,
\end{equation}
under the constraint $\Vert u\Vert _2 =c >0$. A crucial step to
prove that (\ref{1.4}) admits an even ground state $u=u(x)$ which is
decreasing for $x>0$, is (\ref{1.3}), \cite{21}.
\\
In \cite{22}, we showed that the ground state of the
$m\times m$ elliptic system
$$
{\bf (S)}
\left\{
\begin{array}{l}
\Delta u _1 + \lambda _1 u_1 + g_1 (u_1 , \ldots , u_m ) =0
\\
\cdots
\\
\Delta u _m + \lambda _m u_m + g_m (u_1 , \ldots , u_m ) =0
\end{array}
\right.
,
$$
where $g_i = (\partial G / \partial u_i ) $, ($i= 1, \ldots ,m $),
exists when $G$ is supermodular, and satisfies some further conditions.
\\
Notice that supermodularity condition is equivalent to the cooperativity of ${\bf (S)}$.
As proved by W.C. Troy \cite{30}, the latter condition is
also necessary for the existence of solutions of ${\bf (S)}$.
In \cite{10}, G. Carlier viewed the generalized Hardy--Littlewood rearrangement
inequality as an optimal transportation problem. He showed that the
left-hand side of (\ref{1.3}) achieves its maximum (i.e. the cost is minimized)
under the supermodularity assumption. In this paper, we prove, among other things, that (\ref{1.1})
is necessary for the inequality  (\ref{1.3}) to hold.
\\
The generalized Riesz rearrangement inequality is
\begin{eqnarray}
\nonumber
 & & \int \cdots \int F(u_1 (x_1 ), \ldots ,
u_m (x_m )) \prod_{i<j} K_{ij} (d(x_i , x_j )) \, dx_1 \cdots dx_m
\\
\label{1.5}
 & \leq &
\int \cdots \int F(u_1 ^* (x_1 ), \ldots , u_m ^* (x_m ))
\prod_{i<j} K_{ij} (d(x_i , x_j )) \, dx_1 \cdots dx_m ,
\end{eqnarray}
where $d$ denotes  distance and the functions $K_{ij} $ are
decreasing.
\\
It is closely related to the Brunn-Minkowski inequality of convex
geometry. Notice that the integral in (\ref{1.5}) can represent a
physical interaction potential.
\\
The most relevant case is  $K_{ij} (d(x_i , x_j )) =
\omega _{ij} (|x_i -x_j | )$ where $\omega _{ij}$ are nonnegative nonincreasing functions.
(\ref{1.5}) is then equivalent to
\begin{eqnarray}
\nonumber
 & & \int
  \cdots \int F(u_1 (x_1 ), \ldots , u_m (x_m ))
\omega _{ij} (|x_i - x_j |) \, dx_1 \cdots dx_m
\\
\label{1.6}
 & \leq & \int \cdots \int F(u_1 ^* (x_1 ), \ldots , u_m ^* (x_m ))
\omega _{ij} (|x_i - x_j |) \, dx_1 \cdots dx_m .
\end{eqnarray}
(\ref{1.5}) and (\ref{1.6}) have been studied in \cite{28, 29, 3, 4,
5, 1, 26, 14, 27, 12, 31, 7}. All these results use the
supermodularity of $F$. Among the numerous applications of
(\ref{1.6}), let us mention that it was extremely useful to prove
the existence and uniqueness of the minimizing solution of
Choquard's equation \cite{24}. In the special case that $F$ is a
product and $m=2$, (\ref{1.6}) hold on the standard spheres and
hyperbolic spaces \cite{2,3} and it still contains the isoperimetric
inequality as a limit. In this paper we will show (see Proposition
3.2 below) that (\ref{1.1}) is necessary for (\ref{1.6}) to hold.
Notice that , in \cite{7}, cases of equality of (\ref{1.3}) and
(\ref{1.6}) were established under the strict supermodularity
condition, that is to say (\ref{1.1}) with the strict inequality
sign. We will also prove that such a condition is inescapable to
establish cases of equality in (\ref{1.3}) and (\ref{1.6}) (see
Remark 1 below).

\subsection*{2. Notation and preliminaries}
We fix $n\in \mathbb{N} $. Let $\mu $ the Lebesgue measure on $\mathbb{R} ^n$.
If $x\in \mathbb{R} ^n $ and  $r>0$,
let $B_r (x)= \{ y\in \mathbb{R} ^n :\, |y-x|< r\} $.
 For any set $M \subset \mathbb{R} ^n $, let ${\bf 1} _M $
 denote its characteristic function.
By ${\cal M}$ we denote the set of all measurable functions on
$\mathbb{R } ^n $. For a Borel measurable function $F: \mathbb{R} ^m
_+ \to \mathbb{R} $, nonincreasing functions $w_{i j} : \mathbb{R}
_+ \to \mathbb{R} _+ $ and nonnegative $u_i \in {\cal M}$, we study
the following generalized Hardy--Littlewood type functional,
\setcounter{section}{2} \setcounter{equation}{0}
\begin{equation}
\label{2.1}
 I(u_1 , \ldots , u_m ) =
 \int_{\mathbb{R} ^n } F(u_1 (x ), \ldots , u_m (x ))\, dx ,
\end{equation}
and the
generalized Riesz type functional.
\begin{equation}
\label{2.2} J(u_1 , \ldots , u_m ) = \int_{\mathbb{R} ^n } \cdots
\int_{\mathbb{R} ^n}  F(u_1 (x_1 ), \ldots , u_m (x_m )) \omega
_{ij} (| x_i -x_j |) \, dx_1 \cdots dx_m .
\end{equation}
\\
If $f\in {\cal F} $, let $f^*$ denote its {\sl Schwarz
symmetrization} which is the unique lower continuous function which
is radially symmetric, radially nonincreasing and such that $\mu \{
a<f\leq b\} = \mu \{ a<f^* \leq b \}  $ for all numbers $\inf f <a
<b $ (see \cite{1}). Notice that if $f= {\bf 1} _M $, where $M$ is
Lebesgue
 measurable with finite measure,
then $f^* = {\bf 1} _{B_R (0)} $ where $R $ is chosen such that $\mu
(B_R (0))= \mu (M) $ . Accordingly, a function $u\in {\cal  F}$ is
called {\sl Schwarz symmetric } if it is radial and radially
decreasing. We say that it is strictly Schwarz symmetric if it is
radial and strictly radially decreasing.

\section*{ 3. Results}
\setcounter{section}{3}
\setcounter{equation}{0}
{\bf Proposition 3.1}
{\sl ( Necessity of the supermodularity condition in the
generalized Hardy--Littlewood inequality )
\\
Let $F: (\mathbb{R} _+ ) ^m  \to \mathbb{R} $ be a Borel measurable
function which vanishes on hyperplanes, that is,
$F(y_1 , \ldots , y_m ) =0 $ if $y_k =0 $ for some $k\in \{1 , \ldots , m\} $. If
$$
\int_{\mathbb{R} ^n } F(u_1 (x) , \ldots , u_m (x))\, dx \leq
\int_{\mathbb{R} ^n } F(u_1 ^* (x) , \ldots , u_m ^* (x))\, dx,
$$
for any $(u_1 , \ldots , u_m ) \in {\cal F}  ^m $ then $F$ satisfies
(\ref{1.1}).}
\\[0.2cm]
{\sl Proof:} Suppose that (\ref{1.1}) is not true. Then there exist
$ i,j\in \{ 1, \ldots ,m\} $, $i\not=j $, $0\leq a< b$, $0\leq c<d$
and $m-2$ nonnegative numbers $\alpha _l $, ($l\in \{ 1, \ldots , m
\} \setminus \{ i,j\} $),  such that:
\begin{eqnarray}
\label{3.1}
 & & F(\alpha _1 , \ldots , a, \ldots , d, \ldots ,
 \alpha _{m-2} ) + F(\alpha _1 , \ldots , b, \ldots ,c, \ldots , \alpha _{m-2} )
\\
\nonumber
 & > &
F(\alpha _1 , \ldots , b, \ldots , d, \ldots , \alpha _{m-2} ) +
F(\alpha _1 , \ldots , a, \ldots ,c, \ldots , \alpha _{m-2} ) .
\end{eqnarray}
Now let $E$ and $F$ be two measurable sets  such that $E\cap F =
\emptyset $, $\mu (F) = \mu (E) <\infty $, and $u_i  = a {\bf 1} _E
+ b{\bf 1}_F $, $u_j = c {\bf 1}_F  + d {\bf 1}_E  $, and $u_l  =
\alpha _l ({\bf 1}_F  + {\bf 1}_E )$ for $l\in \{ 1, \ldots ,m \}
\setminus \{ i,j \} $. Then $u_i ^*  = b {\bf 1}_{B_r (0)}  + a {\bf
1}_A  $, where $r $ is chosen such that $\mu (B_r (0))= \mu (F) =
\mu (E)$, and $A$ is the annulus $B_{r'} \setminus \overline{B_r } $
and $\mu (A) =\mu (F) = \mu (E)$. It follows that $u_j ^*  = d{\bf
1}_{B_r (0)}  + c {\bf 1}_A $, $u_l ^*  = \alpha _l {\bf 1}_{B_{r'}
(0)} $, and
\begin{eqnarray*}
 & & \int F(u_1 (x) , \ldots , u_n (x)) \, dx
\\
 & = &
[ F(\alpha _1 , \ldots , a, \ldots , d,\ldots , \alpha _{m-2} )
+ F(\alpha _1 , \ldots , b, \ldots , c,\ldots , \alpha _{m-2} )] \mu (E) ,
\\
 & & \int F(u_1 ^*(x) , \ldots , u_n ^*(x)) \, dx
\\
 & = &
[ F(\alpha _1 , \ldots , b, \ldots , d,\ldots , \alpha _{m-2} )
 + F(\alpha _1 , \ldots , a, \ldots , c,\ldots , \alpha _{m-2} )] \mu (E) .
\end{eqnarray*}
Hence we have in view of (\ref{3.1}),
$$
\int F(u_1 (x) , \ldots , u_n (x)) \, dx > \int F(u_1 ^*(x) ,
\ldots , u_n ^*(x)) \, dx,
$$
a contradiction. $\hfill \Box $
\\[0.2cm]
{\bf Remark 1: } Following the same approach we can easily conclude
that the strict supermodularity  assumption, that is, (\ref{1.1})
with strict inequality sign,  is necessary to establish cases of
equality. In \cite{7}, (\ref{1.3}) was proven under the
supermodularity assumption (\ref{1.1}) and an integrability
condition, and cases of equality were established assuming
(\ref{1.1}) with strict inequality sign.
\\[0.2cm]
{\bf Proposition 3.2:} {\sl ( Necessity of supermodularity
assumption in Riesz--type integral )
\\
Let $j:(0,+\infty )  \to \mathbb{R} $ be
a function which is nonincreasing and not identically equal to zero,
and such that
\begin{equation}
\label{limit} \lim_{r\to +\infty } r ^{n-1} j(r) =0 .
\end{equation}
Suppose that $\Psi :[0,+\infty )\times [0,+\infty ) \to \mathbb{R} $
is a Borel measurable function satisfying
\begin{equation}
\label{zeroonhyperplane}
\Psi (u,0)= \Psi (0,v)=0 \quad \forall u\geq 0, \ v\geq 0.
\end{equation}
Finally assume that for all nonnegative functions
 $f,g\in L^{\infty } (\mathbb{R} ^n )  $
with compact support there holds
\begin{equation}
\label{Riesz}
\int_{\mathbb{R}^n}
\int_{\mathbb{R}^n}
\Psi ( f(x), g(y)) j(|x-y|) \, dxdy \leq
\int_{\mathbb{R}^n}
\int_{\mathbb{R}^n}
\Psi ( f^* (x), g^* (y)) j(|x-y|) \, dxdy.
\end{equation}
Then $\Psi $ is supermodular, that is, }
\begin{equation}
\label{supermodular}
\Psi (a+a' , b+b' )-\Psi (a+a' , b) -
\Psi (a, b+b' ) +
\Psi (a,b) \geq 0 \quad \forall a,a',b,b' \in [0,+\infty).
\end{equation}
{\sl Proof: } In view of (\ref{limit}),  and since $j$ is
non-increasing and nontrivial, $j$ is non-negative, too, and  there
exist  positive numbers $\varepsilon $ and $t_0 $ with  $t_0 >2
\varepsilon $ and such that
\begin{equation}
\label{strict} j(t)<j(s) \quad \mbox{ whenever $0\leq s\leq
\varepsilon $ and $|t-t_0 | \leq \varepsilon $.}
\end{equation}
We fix $z\in \mathbb{R}  ^n $ with $|z|= t_0 $, and $\rho >t_0 + \varepsilon $.
\\
Now
let $a,a', b,b' \in [0, +\infty )$.
Then, if  $f= a{\bf 1}_{B_R (0)} + a' {\bf 1}_{B_{\varepsilon } (0)}
$, and
$g= b{\bf 1}_{B_{\rho} (0)} + b' {\bf 1}_{B_{\varepsilon } (z)} $,
where $R>\varepsilon $, we have that  $f=f^* $,
and $g^* = b{\bf 1}_{B_{\rho} (0)} + b' {\bf 1}_{B_{\varepsilon } (0)} $. Using
(\ref{zeroonhyperplane}), a short computation shows that
(\ref{Riesz}) is equivalent to
\begin{eqnarray}
\nonumber
0  & \leq & \left\{
\Psi (a+a' , b+b' )-\Psi (a+a' , b) - \Psi (a, b+b' ) +
\Psi (a,b) \right\} \cdot
\\
\label{ineq1}
 & &
\cdot \int_{B_{\varepsilon } (0)} \left[ \int_{B_{\varepsilon } (0)}
j(|x-y|)\, dy - \int_{B_{\varepsilon } (z)} j(|x-y|)\, dy \right] dx
\\
\nonumber
 & + & \left( \Psi (a,b+b' )-\Psi (a,b) \right)
\int_{B_R (0)} \left[ \int_{B_{\varepsilon } (0)} j(|x-y|)\, dy -
\int_{B_{\varepsilon } (z)} j(|x-y|)\, dy \right] dx .
\end{eqnarray}
Set
$$
H(x):=
 \int_{B_{\varepsilon } (0)} j(|x-y|)\, dy , \quad \, x\in \mathbb{R} ^n .
$$
Since $j$ is nonnegative and  nonincreasing,
$H$ is nonnegative, radial and radially nonincreasing,
that is we can write
$H(x) =: h(|x|) $, with $h$ nonincreasing.
Since $j$ satisfies (\ref{limit}),
we also have $\lim _{t\to \infty } t^{n-1} h(t) =0$.
Setting
$$
I(t) := \int_{B_t (0)} H(x)\, dx - \int _{B_t (z) } H(x)\, dx , \quad (t>0 ),
$$
it follows that $I$ is nonnegative, and moreover, (\ref{strict})
implies that
\begin{equation}
\label{strict2}
I(\varepsilon ) >0.
\end{equation}
Finally, if $t>t_0 $, then
 $$
 I(t)  \leq  \int_{B_{t + t_0 } (0) \setminus B_{t-t_0 } (0)} H(x) \, dx
   =  n \omega _n \int_{t-t_0 } ^{t+ t_0 } s^{n-1} h(s) \, ds,
  $$
  where $\omega _n $ denotes the volume of the unit ball in $\mathbb{R} ^n $.
Hence
\begin{equation}
\label{limI}
\lim_{t\to \infty } I(t) =0.
\end{equation}
Now inequality (\ref{ineq1}) becomes
\begin{eqnarray*}
0  & \leq & \left\{
\Psi (a+a' , b+b' )-\Psi (a+a' , b) - \Psi (a, b+b' ) +
\Psi (a,b) \right\} I (\varepsilon )
\\
\label{ineq2}
 & + & \left( \Psi (a,b+b' )-\Psi (a,b) \right) I(R) .
\end{eqnarray*}
Sending $R\to +\infty $, we obtain
\begin{equation}
\label{ineq3} 0  \leq \left\{ \Psi (a+a' , b+b' )-\Psi (a+a' , b) -
\Psi (a, b+b' ) + \Psi (a,b) \right\} I (\varepsilon ),
\end{equation}
and (\ref{supermodular}) follows from (\ref{strict2}).
$\hfill \Box $
\\[0.2cm]
{\bf Remark 2: } Proposition 3.2 tells us in particular that $\Psi $
is supermodular if (\ref{Riesz}) holds for $j = {\bf 1} _{B_R (0)}
$, with some $R>0 $.
\\[0.2cm]
 {\bf Proposition 3.3.} {\sl Let $\Psi :
[0,+\infty ) \times [0,+\infty ) $ be nontrivial, and satisfies
(\ref{supermodular}) and (\ref{zeroonhyperplane}), and let $h\in C
(\mathbb{R} ^n ) $. Assume that there holds
\begin{equation}
\label{riesz2} \int_{\mathbb{R} ^n} \int_{\mathbb{R} ^n } \Psi
(f(x), g(y)) h(x-y ) \, dxdy \leq \int_{\mathbb{R} ^n}
\int_{\mathbb{R} ^n } \Psi (f ^* (x), g ^* (y)) h(x-y ) \, dxdy
\end{equation}
for all bounded nonnegative functions $f$ and $g$ with compact support .
Then $h$ is radial and radially nonincreasing, that is, it can be
written as $j(|x|) = h(x)$, ($x\in \mathbb{R} ^n \setminus \{0 \}  $), where $j:(0,
+\infty ) \to \mathbb{R} $ is a nonincreasing function.  }
\\[0.2cm]
{\sl Proof :} By choosing $a=b=0 $ in (\ref{supermodular}) and
taking into account (\ref{zeroonhyperplane}), we have that $\Psi $
is nonnegative, and since $\Psi $ is nontrivial, there are numbers
$a>0 $, $b>0 $ such that $\Psi (a,b) >0 $. Let $z_1 , z_2 \in
\mathbb{R } ^n  $  with $0< |z_1 |< |z_2 | $, and set $R:= (|z_1 | +
|z_2 | )/2 $.  Let $f= a( {\bf 1}_{ B_R (0)} + {\bf 1 }
_{B_{\varepsilon } (z_2 ) } - {\bf 1}_{B_{\varepsilon } (z_1 ) } )$,
and $g= b {\bf 1}_{B_{\varepsilon } (0) }$, where $0< \varepsilon
<(|z_2 |-|z_1 |)/2 $. Then $g=g^* $ and $f^* = a {\bf 1} _{B_R (0) }
$. Hence (\ref{riesz2}) gives
$$
0 \leq \Psi (a,b)
\left\{ \int_{B_{\varepsilon } (z_1 )}
\int_{B_{\varepsilon } (0)}
h(x-y)\, dxdy - \int_{B_{\varepsilon } (z_2 )}
\int_{B_{\varepsilon } (0)} h(x-y)\, dxdy \right\} .
$$
Sending $\varepsilon \to 0 $
and taking into account
$\Psi (a,b) >0$ we find that $h(z_1 ) \geq h(z_2 ) $.
Since $h$ is continuous,
this implies $h(x)\geq h(y) $ iff
$|x|\leq |y|$.
Hence $h$ is radial and radially nonincreasing as claimed.
$\hfill \Box$
\\[0.2cm]
{\bf Acknowledgement:} 
The  author  is  very grateful to F. Brock for his kind invitation to visit  AUB in June 2008  and June 2009, and for his  precious  help.

\end{document}